\documentstyle[12pt]{article}

  \newcommand{\const}{\rm const}
  \newcommand{\Var}{\rm Var}

\begin{document}

 \begin{center}

 \ {\bf Exponential confidence interval  based on the recursive }\par

\vspace{4mm}

{\bf Wolverton - Wagner  density estimation.} \par

\vspace{5mm}

{\bf  M.R.Formica, E.Ostrovsky, and L.Sirota. } \par

 \vspace{5mm}

 \end{center}

 \ Universit\`{a} degli Studi di Napoli Parthenope, via Generale Parisi 13, Palazzo Pacanowsky, 80132,
Napoli, Italy. \\

e-mail: mara.formica@uniparthenope.it \\

\vspace{3mm}

Department of Mathematics and Statistics, Bar-Ilan University,
59200, Ramat Gan, Israel. \\

\vspace{4mm}

e-mail: eugostrovsky@list.ru\\
Department of Mathematics and Statistics, Bar-Ilan University,\\
59200, Ramat Gan, Israel.

\vspace{4mm}

e-mail: sirota3@bezeqint.net \\

\vspace{5mm}

\begin{center}

 \ {\bf Abstract.}

 \end{center}

 \ We  derive the exponential  non improvable  Grand Lebesgue Space norm decreasing estimations for tail of distribution for
  exact normed deviation for the famous recursive  Wolverton - Wagner  multivariate  statistical density estimation. \par
 \ We consider pointwise  as well as  Lebesgue - Riesz norm error of statistical density of measurement. \par

 \vspace{5mm}

  \par

\vspace{5mm}

\section{Statement of problem. Notations and definitions. Previous results.}

\vspace{5mm}

 \ Let  $ \ (\Omega,  M, {\bf P})  \ $  be probability space with expectation $ \ {\bf E} \ $ and variance $ \ \Var. \ $
 \ Let also  $ \ \{\xi_k \}, \ k = 1,2,\ldots,n \ $ be a sequence of independent, identical distributed (i, i.d.) random {\it vectors } (r.v.)
 taking the values in the ordinary Euclidean space $ \ R^d, \ d = 1,2,\ldots \ $ and having certain {\it non - known} density of a distribution
 $ \ f = f(x), \ x \in R^d. \ $ C.Wolverton and T.J.Wagner in \cite{Wolverton Wagner} offered the following famous statistical estimation $ \ f^{WW}_n(x) = f_n(x) \ $
 for $ \ f(\cdot). \ $ \par
  \ Let $ \ \{h_k\}, \ k = 1,2,\ldots \ $ be some positive sequence of real numbers such that $ \ \lim_{k \to \infty} h_k = 0. \ $ Let also $ \ K = K(x), \ x \in R^d \ $
 be certain {\it kernel},  i.e. measurable even function  for which

\begin{equation} \label{Kernel}
\int_{R^d} K(x) \ dx = 1.
\end{equation}
 \ Then by definition

\begin{equation} \label{WW est}
f^{WW}_n(x) = f_n(x) \stackrel{def}{=} \frac{1}{n} \sum_{k=1}^n \frac{1}{h_k^d} \ K \left(\ \frac{x - \xi_k}{h_k}    \ \right).
\end{equation}

 \ Recall that the classical kernel, or Parzen - Rosenblatt's  estimate $ \ f^{PR}_n(x) \ $ has a form

\begin{equation} \label{OPR est}
f^{PR}_n(x) := \frac{1}{nh^d_n} \sum_{k=1}^n K \left( \   \frac{x - \xi_k}{h_n}  \  \right),
\end{equation}
 see \cite{Parzen}, \cite{Rosenblatt}. \par
 \ Note that the Wolverton - Wagner estimate obeys a very important recursion property:

$$
f^{WW}_n(x) = \frac{n-1}{n} \ f^{WW}_{n-1}(x) + \frac{1}{n h_n^d} \ K \left(  \  \frac{x - \xi_n}{h_n}  \ \right).
$$

\vspace{4mm}

 \ "The recurrent definition of probability density estimates" $ \ f_n^{WW}(x) \ $
 has two obvious advantages: 1) there is no need to memorize data, i.e. if the estimate $ \ f^{WW}_{n-1}(x) \ $
 is known, then $ \ f_n^{WW}(x) \ $ can be calculated by means of the last
observation $ \ \xi_n \ $   only, without using the sampling $ \ \xi_1, \ \xi_2, \ldots, \xi_{n-1};  \ $ 2) the asymptotic
dispersion of the estimate $ \ f^{WW}_n(x) \ $  does not exceed the dispersion of the estimate" $ \ f^{PR}(x), \ $  see
\cite{Nadaraya}. \par

\vspace{4mm}

 \ {\bf Our aim in this report is to deduce the  exact {\it exponential decreasing}  estimate for the tail of deviation probability}

\begin{equation} \label{statement tail}
{\bf P}_n^{WW}(u) \stackrel{def}{=}  \sup_{x \in R^d}  {\bf P} (B_n |f^{WW}_n(x) - f(x)| > u \ ), \ u \ge 1,
\end{equation}
i.e. under {\it exact} optimal deterministic numerical sequence $ \ B_n, \ $ such that $ \ \lim_{n \to \infty} B_n = \infty. \ $\par

\vspace{3mm}

 \ For the Parzen - Rosenblatt estimate $ \ f_n^{PR}(x) \ $ these estimates was obtained e.g. in \cite{Ostrovsky0}, chapter 5,
 sections 5.2 - 5.6. \par
 \ We will use some facts from the theory of the so - called Grand Lebesgue Spaces (GLS), devoted in particular the Banach spaces
 of random variables having exponential decreasing tails of distributions, see e.g.
\cite{anatrielloformicaricmat2016} , \cite{Buldygin-Mushtary-Ostrovsky-Pushalsky},  \cite{caponeformicagiovanonlanal2013},
\cite{Iwaniec1}, \cite{Iwaniec2}, \cite{Kozachenko},  \cite{Ostrovsky0} etc. \par

 \ Note that the distribution of the normed deviation  $ \ B_n(f^{WW}_n - f) \ $ in different Lebesgue - Riesz spaces
 $ \ L_p(R^d \otimes \Omega ): \ $

$$
\Delta_{n,p} := {\bf E} \int_{R^d} B_n^p | \  f^{WW}_n(x) - f(x)  \ |^p \ dx
$$
 was investigated in many works, see e.g. \cite{Devroe}, \cite{Nadaraya}, \cite{Parzen}, \cite{Rosenblatt},
 \cite{Slaoui},  \cite{Tsybakov1}, \cite{Tsybakov2}, \cite{Wegman}, \cite{Wolverton Wagner} etc. The optimal choose of
$ \ \{h_k\} \ $ and the kernel $ \ K(x) \ $ are devoted the following works
\cite{Compte},   \cite{Goldenshluger Lepski},  \cite{Khardani}, \cite{Lerasle}. The case when the r.v. - s. are dependent
is investigated in  \cite{Khardani}, \cite{Simone}. \par

\vspace{4mm}

\begin{center}

   \ {\sc Let us  reproduce some used  for us notations and conditions from this theory.} \par

\vspace{3mm}

\end{center}

 \ Let $ \ (\beta, L)  \ $ be certain positive numbers. Denote by $ \ l= l(\beta) = [\beta] \ $ an integer part for $ \ \beta, \ $ i.e.
 a maximal integer number less than $ \ \beta: \ $

$$
 l(\beta) = [\beta] = \max \{ \ j = 0,1,2,\ldots: j \le \beta \ \}.
$$
 \ For instance, $ \ l(0. 3) = 0, \ l(\pi) = 3. \ $  Correspondingly, the fractional part $ \ \{\beta\} \ $  for $ \ \beta \ $ is equal
 $ \ \{\beta\} := \beta - [\beta]. \ $ \par
 \ Introduce as ordinary the functional class $ \ \Sigma(\beta,L) \ $ as follows

\begin{equation} \label{sigma beta L}
\Sigma(\beta,L) = \{ \ f: R^d \to R; \ \forall m = \vec{m}: \ |m| \le [\beta] \ \Rightarrow \frac{\partial^m f}{\partial x^m} \in H(\{\beta\},L\ ), \ \}
\end{equation}
where $ \ H(\alpha,L) \ $ denotes the H\"older class of the functions

$$
H(\alpha,L) = \{g: R^d \to R, \ |g(x) - g(y)| \le L \cdot |x-y|^{\alpha} \ \}, \ \alpha \in (0,1]
$$
 \ As usually
$$
  |z| = \sqrt{(z,z)}, \ m = \vec{m} = \{m_1,m_2, \ldots,m_d \}, \ |m| = \sum_{i=1}^d m_i.
$$
 \ In the case when $ \ \beta \ $ is integer number, the derivative in  (\ref{sigma beta L}) is assumed  to be continuous and bounded:

$$
 \forall m = \vec{m}: \ |m| = \beta = [\beta] \ \Rightarrow \sup_{x \in R^d} \left| \frac{\partial^{\vec{m} f}}{\partial x^{\vec{m}}} \ \right|  \le L.
$$

\vspace{4mm}

\ {\it We suppose henceforth  that the density function  belongs to some set   } $ \ \Sigma(\beta,L) \ $
for some non - trivial value $ \ \beta \in (0,\infty): \ $

\begin{equation} \label{Sigma beta}
f(\cdot) \in \Sigma(\beta) \stackrel{def}{=} \cup_{L \in (0,\infty)} \Sigma(\beta,L), \ \beta > 0.
\end{equation}

\vspace{4mm}

 \ As for the kernel $ \ K. \ $ We impose on $ \ K \ $ the following conditions

\begin{equation} \label{K cond 1}
K(-x) = K(x); \ \int_{R^d} K(x) \ dx = 1; \ \int_{R^d} K^2(x) dx < \infty;
\end{equation}

\begin{equation} \label{K cond 2}
K(\cdot) \in C(R^d), \ \int_{R^d} |K(x)| dx < \infty.
\end{equation}

 \ The following conditions may be named as  {\it conditions of orthogonality:}

\begin{equation} \label{cond orthog}
\forall \vec{m}: \ |\vec{m}| \le [\beta] \ \Rightarrow \int_{R^d} \vec{x}^{\vec{m}} \ V(x) \ dx = 0.
\end{equation}

 \ The  last conditions (\ref{cond orthog})
  may be used only for the investigation of {\it bias} $ \ \delta_n(x) \ $ of these statistics

\begin{equation} \label{bias}
\delta_n =  \delta_n(x) \stackrel{def}{=} {\bf E} f_n^{WW}(x) - f(x).
\end{equation}

 \ In detail, as long as $ \ f \in \Sigma(\beta) \ $ and by virtue  of  (\ref{cond orthog})

$$
\delta^{(k)} \stackrel{def}{=} h_k^{-d} \int_{R^d} K \left( \ \frac{x - y}{h_k} \ \right)  \ f(y) \ d y - f(x) \sim h_k^{\beta}, \ h_k \to 0+,
$$
therefore

\begin{equation} \label{delta n est}
|\delta_n| \sim C_1(\beta,L) \ n^{-1} \left[ \ \sum_{k=1}^n h_k^{\beta} \ \right],
\end{equation}

see \cite{Nadaraya},   \cite{Devroe},  \cite{Khardani}. \par

\vspace{5mm}

\section{Main result.}

\vspace{5mm}

 \ Let us investigate now the {\it  Variance} of the considered Wolver- Wagner  $ \ f^{WW}_n(x) \ $  statistic, of  course under
 formulated above restrictions. We have

$$
\Var \left\{ \ f^{WW}_n(x) \ \right\} = \frac{1}{n^2} \ \sum_{k=1}^n  \ h_k^{-d} \  \Var \left\{ \ K \left( \ \frac{x - \xi_k}{h_k}   \ \right) \ \right\} \asymp
\frac{1}{n^2} \sum_{k=1}^n \frac{1}{h_k^d}.
$$

 \ Let's form the classical target functional

$$
Z = Z_n(h_1,h_2, \ldots, h_n)  \stackrel{def}{=} {\bf E} \left[ \ f_n^{WW}(x) - f(x) \ \right]^2;
$$
then

\begin{equation} \label{Z fun}
Z_n(h_1,h_2, \ldots, h_n)  \asymp  \frac{1}{n^2} \   \left\{ \ \sum_{k=1}^n \frac{1}{h_k^d} + \left[ \ \sum_{k=1}^n h_k^{\beta}  \ \right]^2 \ \right\}.
\end{equation}

 \ The (asymptotic) minimal value of the functional $ \ Z_n(h_1,h_2, \ldots, h_n) \ $ relative the variables $ \ \{h_k\} \ $ subject to our limitations
is attained on the values

\begin{equation} \label{h k}
h_k \sim  C_2(\beta,d,L) \ k^{-1/(2 \beta + d)}
\end{equation}
 and wherein

\begin{equation} \label{min Z}
\min Z_n = n^{- 2 \beta/(2 \beta + d)}.
\end{equation}

 \vspace{3mm}

  \ So, the speed of convergence $ \ f_n^{WW}(x) \to f(x) \ $ as $ \ n \to \infty \ $ is equal to $ \ n^{- \beta/(2 \beta + d)}:  \ $

\begin{equation} \label{speed}
\left[ \  {\bf E} (f_n^{WW}(x) - f(x))^2 \  \right]^{1/2}  \asymp n^{-\beta/(2 \beta + d)},
\end{equation}
alike ones  for the Parzen - Rosenblatt estimates. \par

 \ On the other words, the value $ \ B_n \ $ in (\ref{statement tail}) must be choosed as follows:

\begin{equation} \label{Bn opt}
  B_n = n^{\beta/(2 \beta + d)}.
\end{equation}

\vspace{3mm}

\ Note that the one - dimensional case $ \ d = 1 \ $ was considered in \cite{Compte}. \par

\vspace{4mm}

 \ {\it  We suppose henceforth that the values  } $ \ \{h_k\}, \ B_n \ $ {\it are choosed optimally in accordance with }
(\ref{h k}) {\it and } (\ref{Bn opt}). \par

\vspace{4mm}

 \ Define the following tail probability

\begin{equation} \label{Qn tail}
{\bf Q}_n^{WW}(u) \stackrel{def}{=}  \sup_{x \in R^d}  {\bf P} (B_n |f^{WW}_n(x) - {\bf E}f^{WW}_n(x)| > u \ ), \ u \ge 1.
\end{equation}

 \ Note that the following value is bounded:

$$
\sup_x \ \sup_n  B_n |f(x) - {\bf E}f^{WW}_n(x)| = C_3 = C_3(d, \beta,L) < \infty.
$$
 \ Therefore

$$
{\bf P}_n^{WW}(u) \le  {\bf Q}_n^{WW}(u - C_3).
$$
 \ Evidently, the r.v. $ \ f^{WW}_n(x) - {\bf E}f^{WW}_n(x) \ $ is centered (mean zero). \par

 \vspace{4mm}

 \ {\bf Theorem 2.1.} We propose under formulated above conditions

\begin{equation} \label{main}
  \sup_n  \ {\bf Q}_n^{WW}(u) \le 2 \exp \left[ \ - C_4(d,\beta,L) \ u^{ \frac{2\beta + d}{ \beta + d}  } \  \right], \ u \ge 1.
\end{equation}

\vspace{4mm}

 \ {\bf Proof.} First of all we need for applying  the theory of Grand Lebesgue Spaces (GLS) the estimate
 of an exponential moment

\begin{equation} \label{exp mom}
E_n[Q,\lambda, \beta] \stackrel{def}{=} {\bf E} \exp \left[ \ \lambda B_n (f^{WW}_n(x) - {\bf E}f^{WW }_n(x)) \ \right], \ \lambda \in R.
\end{equation}

\ Denote for this purpose

$$
\Theta_n = B_n (f^{WW}_n(x) - {\bf E}f^{WW }_n(x)) = n^{- \frac{ \beta + d}{2 \beta +d} } \ \sum_{k=1}^n h_k^{-d} K^o \left( \ \frac{x - \xi_k}{h_k}  \ \right) =
$$

$$
\sum_{k=1}^n \theta_{k,n}, \ \theta_{k,n} :=   n^{- \frac{ \beta + d}{2 \beta +d} } \  h_k^{-d} K^o \left( \ \frac{x - \xi_k}{h_k}  \ \right),
$$
where as ordinary for arbitrary r.v. $ \ \eta \ \Rightarrow \ \eta^o \stackrel{def}{=} \eta - {\bf E} \eta. \  $ We have

$$
E_n[Q,\lambda, \beta] = {\bf E}e^{\lambda \Theta_n} = \prod_{k=1}^n {\bf E} e^{\lambda \theta_{k,n}}.
$$

 \ Let us consider two possibilities.

$$
 \ {\bf A.} \hspace{5mm} 0 <  \lambda \le n^{- \frac{\beta + d}{2 \beta + d} } \ h_n^{-d} /(2 \sup_x K(x)) \le 1,
$$
or equally

$$
\lambda  \in \left( \ 0, \ C_5 n^{\frac{\beta}{2 \beta + d} } \ \right).
$$

 \ We use the following elementary inequality

$$
y \in (0,1) \ \Rightarrow e^y \le 1 + y + y^2.
$$

 \ Therefore

$$
E_{k,n}(\lambda) := {\bf E} \exp \left( \ \lambda \theta_{k,n} \ \right) \le 1  + C \lambda^2 \Var (\theta_{k,n}) \le
$$

$$
 \exp \left( \ C  \lambda^2 \Var (\theta_{k,n})  \ \right);
$$

$$
E_n[Q,\lambda, \beta]\le \exp \left( \ C \lambda^2 \sum_{k=1}^n \Var (\theta_{k,n}) \ \right) \le \exp (C \lambda^2).
$$

\vspace{3mm}

 \ {\bf B.} Let us investigate an opposite possibility

$$
\lambda  \ge \ C_5 \ n^{\frac{\beta}{2 \beta + d} }.
$$

 \ But then

$$
n \le C \ \lambda^{ \  \frac{2\beta + d}{\beta} \ },
$$
and we deduce

$$
\lambda \Theta_n = \lambda \ n^{  \ - \frac{\beta + d}{2 \beta + d}  \ } \ \sum_{k=1}^n h_k^{-d} \ K^o \left( \ \frac{x - \xi_k}{h_k} \ \right),
$$
and following

$$
\lambda |\Theta_n|  \le C \lambda n^{ \ \frac{\beta + d}{2 \beta + d} \ } \le C \lambda^{ \ \frac{2 \beta + d}{\beta}  \ },
$$

$$
{\bf E} \exp ( \lambda \Theta_n) \le \exp \left( \  C \lambda^{ \ \frac{2 \beta + d}{\beta}  \ }   \ \right).
$$

\vspace{3mm}

 \ The case when $ \ \lambda < 0 \ $ is considered quite analogously. \par

\vspace{4mm}

 \ Denote $ \ m = m(n) = n^{\beta/(2 \beta + d)}. \ $  Let us introduce the following function

$$
\phi(\lambda)= \phi_m(\lambda) = \phi_{\beta,d,n}(\lambda) = \lambda^2 I( |\lambda| \le m  ) +
|\lambda|^{(2 \beta + d)/\beta} I( |\lambda| > m  ),
$$
where as ordinary $ \ I(A) \ $ denotes the indicator function of the set $ \ A. \ $ \par
 \ We obtained actually

\begin{equation} \label{exp mom}
\ln {\bf E} \ \left( \lambda \ \Theta_n \right) \le \phi_m(C_6 \lambda), \lambda \in R.
\end{equation}

 \ It follows from the theory of Grand Lebesgue Spaces (GLS), see e.g.  \cite{Buldygin-Mushtary-Ostrovsky-Pushalsky},
\cite{caponeformicagiovanonlanal2013}, \cite{Formica etc}, \cite{Kozachenko} that

$$
{\bf P}(|\Theta_n| > C u) \le \exp \left( \ -  \phi^*_m(u)   \ \right), \ u \ge 0,
$$
modified Tchernov's  inequality. Here as usually $ \ \phi^*(\cdot) \ $ denotes the classical Young - Fenchel transform

$$
\phi^*(u)  \stackrel{def}{=} \sup_{\lambda \in R} (\lambda u - \phi(\lambda)).
$$

 \ We deduce after simple calculations

\begin{equation} \label{complete estimation begin}
{\bf P}(|\Theta_n| > C u) \le \exp \left( \ -  u^2 \ \right), \ u \in \left( \ 0, n^{\beta/( 2 \beta + d)}  \ \right);
\end{equation}

\begin{equation} \label{complete estimation end}
{\bf P}(|\Theta_n| > C u) \le \exp \left( \ -  u^{(2\beta + d)/(\beta + d)} \ \right), \ u \ge  n^{\beta/( 2 \beta + d)}.
\end{equation}

\vspace{3mm}

 \ The announced result (\ref{main}) follows immediately from (\ref{complete estimation begin}) and (\ref{complete estimation end});
 it is easily to verify that obtained estimate for  $ \  {\bf P}(|\Theta_n| > C u)  \ $   reaches its maximum relative the variable $ \ n \ $
 only for the value $ \ n = 1. \ $ \par

\vspace{3mm}

 \ {\bf Remark 2.1.} \ Note that the inequality  of the form (\ref{main}) of  Theorem 2.1 is true also for the classical Parzen - Rosenblatt estimation,
see \cite{Ostrovsky0}, chapter 5, sections 1 - 2. \par

\vspace{3mm}

 \ {\bf Remark 2.2.}  \ It is known, see \cite{Ostrovsky0}, chapter 5, section 3 that the result  (\ref{main}) is essentially non - improvable. Indeed,
 there holds the following {\it lower} estimate under our conditions for {\it arbitrary} density statistics  $ \ \hat{f}: \ $

$$
  \sup_n \ \sup_x  \ {\bf P} \left( \ B_n |\hat{f}(x) - f(x)| > u \ \right) \ge
$$

\begin{equation} \label{lower}
 \ge 2 \exp \left[ \ - C_{14}(d,\beta,L) \ u^{ \frac{2\beta + d}{ \beta + d}  } \  \right], \ u \ge 1.
\end{equation}

\vspace{5mm}

\section{Error estimate in Lebesgue - Riesz norms. }

\vspace{5mm}

 \ Let $ \ \mu \ $ be arbitrary Borelian {\it finite: } $ \ \mu(R^d) = 1 \ $ measure on the whole space  $ \ X := R^d. \ $
For instance,

$$
\mu(A) = \frac{\nu(A \cap D)}{\nu(D)},
$$
where $ \ \nu \ $ is ordinary Lebesgue measure and $ \ D \ $  is fixed measurable non - trivial set: $ \ 0 < \nu(D) < \infty \ $
or $ \ \mu = \delta_{x_0} - \ $ unit delta Dirac measure concentrated at the point $ \ x_0 \in X. \ $ \par

 \ Introduce as ordinary the classical Lebesgue - Riesz space $ \ L_p = L_p(R^d,\mu) \ $ as a set of all the (measurable)
functions $ \ g: R^d \to R \ $ having a finite norm

$$
||g||_p := \left[  \ \int_{R^d} |g(x)|^p \ \mu(dx)  \ \right]^{1/p}, \ p \in [1,\infty).
$$

 \ {\it  We intent in this section to evaluate the error estimation of the Wolverton - Wagner statistics in the $ \ L_p \ $ norm:}

 \vspace{3mm}

\begin{equation} \label{Lp statement}
R_{n,p}(u) \stackrel{def}{=}  {\bf P} ( \ B_n ||f^{WW}_n - f||_p > u), \ u \ge 1.
\end{equation}

\vspace{3mm}

 \ Note that case $ \ p = 2 \ $ (Hilbert space) was considered in many works,  e.g.  \cite{Compte}, \cite{Devroe},
\cite{Nadaraya} at all. The $ \ L_1 \ $  approach was investigated in the monograph \cite{Devroe Giorfi}. \par

\vspace{4mm}

\ {\bf Theorem 3.1.} We propose again under formulated above conditions

\begin{equation} \label{main Lp}
  \sup_n \ R_{n,p}(u)  \le  \exp \left[ \ - C_8(d,\beta,L,p) \ (u- C_3)^{ \frac{2\beta + d}{ \beta + d}  } \  \right], \ u \ge C_3.
\end{equation}

\vspace{4mm}

 \ {\bf Proof.} \ We need for this purpose to apply the theory of the so - called {\it mixed } Lebesgue - Riesz spaces, e.g.
 \cite{Benedek  Panzone}, \cite{Besov Ilin Nikolskii}. Indeed, let us introduce the following two mixed Lebesgue - Riesz spaces
 containing on all the bi-measurable  numerical valued random processes (fields)
 $ \ \eta = \eta( x,\omega), \ x \in R^d, \ \omega \in \Omega \ $ having a finite norm

\begin{equation} \label{pX rOmega beg}
||\eta||_{p,X,r,\Omega} \stackrel{def}{=} || \   || \ \eta \ ||_{p,X} \ ||_{r,\Omega} =
\end{equation}

\begin{equation} \label{pX rOmega end}
\left\{ \ {\bf E} \left[ \ \int_{X} |\eta(x,\omega)|^p \ \mu(dx) \ \right]^{r/p}  \ \right\}^{1/r},  \ X = R^d, \ p,r \ge 1.
\end{equation}
and correspondingly

\vspace{3mm}

\begin{equation} \label{rOmega pX beg}
||\eta||_{r,\Omega, p,X} \stackrel{def}{=} || \   || \ \eta \ ||_{r,\Omega} \ ||_{p, X} =
\end{equation}

\begin{equation} \label{rOmega pX end}
\left\{ \  \int_X \left[ \  {\bf E} |\eta(x)|^r  \ \right]^{p/r} \ \mu(dx) \ \right\}^{1/p},  \ X = R^d, \ p,r \ge 1.
\end{equation}

 \ Evidently, in general case $ \ ||\eta||_{p,X,r,\Omega} \ne ||\eta||_{r,\Omega, p,X}, \ $ but always

\begin{equation}  \label{key inequality}
||\eta||_{p,X,pr,\Omega} \le ||\eta||_{pr,\Omega,p,X}.
\end{equation}

\vspace{3mm}

 \ It  follows from the theory of Grand Lebesgue Spaces, \cite{Kozachenko}, \cite{Ostrovsky0}, chapter 1, section 1.5 that if
 the r.v. $ \ \zeta \ $ satisfies the inequality

$$
{\bf P}(|\zeta| > u) \le \exp \left( \ - C u^{(2 \beta+ d)/(\beta + d)} \ \right),
$$
 then

$$
\sup_{r \ge 1}  \left[ \  r^{- \beta/(2 \beta + d)} ||\zeta||_{r,\Omega} \ \right] < \infty,
$$
 and inverse proposition is also true. Therefore   it follows from the estimation   (\ref{complete estimation end})
 that uniformly relative the parameter $ \ n \ $

$$
||\Theta_n||_{r,\Omega} \le C_9(\beta,L,d) r^{\beta/(2 \beta + d)}, \ r \ge 1..
$$

 \ We obtain from the relation (\ref{key inequality}) denoting $ \ \kappa = ||\Theta_n||_{p,X}: \ $

$$
||\kappa||_{pr, \Omega} \le C_{10} (pr)^{ \ \beta/(2 \beta + d) \ },
$$
or equally taking into account the boundedness of the measure $ \ \mu \ $ and applying Lyapunov's inequality

$$
||\kappa||_s \le C_{11}(\beta,d,L;p) \ s^{\beta/(2 \beta + d)}, \ s \ge 1;
$$
which is completely  equal the the  assertion of theorem 3.1. \par

\vspace{4mm}

 \ {\bf Corollary 3.1.} It follows from theorem 3.1

\begin{equation} \label{estim}
{\bf P}(||f^{WW}_n - f||_p > C (v + C_3)) \le \Delta_n(v),
\end{equation}
where

\begin{equation} \label{Delta n}
\Delta_n(v) := \exp \left( \ - n^{\beta/(\beta + d)} \times v^{(2 \beta + d)/(\beta + d)} \ \right), \ v = \const \ge 1.
\end{equation}
 \ Since

$$
\forall v > 1 \ \Rightarrow \sum_{n=1}^{\infty} \Delta_n(v) < \infty,
$$
we conclude that for all the values $ \ p \in [1,\infty) \ $ the Wolverton - Wagner's statistics converges in the $ \ L_p \ $
 norm with probability one:

\begin{equation} \label{converg}
{\bf P} ( \ ||f^{WW} - f||_p \to 0 \ )  = 1.
\end{equation}

\vspace{4mm}

 \ {\bf Corollary 3.2.} Moreover, define for each  the value$ \ p \in [1, \infty) \ $ the variables

$$
\tau_n \stackrel{def}{=} ||f^{WW}_n - f||_p, \
\tau \stackrel{def}{=} \sup_n ||f_n^{WW} - f||_p = \sup_n \tau_n.
$$
 \ As long as

$$
{\bf P}(\tau > v) \le \sum_{n=1}^{\infty} {\bf P}(\tau_n > v), \ v \ge 1,
$$
we get after some calculations

$$
{\bf P}(\tau > C_{12}( v + C_3) ) \le \sum_{n=1}^{\infty}  \exp \left( \ - v^{ (2 \beta + d)/(\beta + d) } \cdot n^{\beta/(\beta + d)} \ \right) \le
$$

$$
 C_{14} v^{- (2\beta + d)/\beta }, \ v \ge 1.
$$

 \vspace{5mm}

\section{Concluding remarks.}

 \vspace{5mm}

 \hspace{3mm} {\bf A.} It is interest by our opinion to deduce the optimal density estimation as well as
 confidence region  in the uniform norm $ \ L_{\infty}: \ $

 $$
 ||g||_{\infty} := \sup_{x \in R^d} |g(x)|.
 $$
 \ Perhaps, one can set for this purpose

$$
h_k := \frac{(\ln k)^{\gamma}}{k^{1/(2 \beta + d) }}.
$$

\hspace{3mm}

 \ {\bf B.} Offered here method may be generalized on the so called regression problem, i.e. when

$$
\eta_i = f(x_i) + \epsilon_i, \ i = 1,2,\ldots,n.
$$

\hspace{3mm}

 \ {\bf C.}  For the practical using it may be interest to investigate the weight approximate for density, e.g.

$$
\Gamma_{n,w}[f_n^{WW}, f] := \sup_n \sup_x [ \ w(x) \ B_n \ ||f_n^{WW}(\cdot) - f(\cdot)|| \ ],
$$
where $ \ w = w(x) \ $ is certain weight, i.e. non - negative numerical valued measurable function. \par

\vspace{6mm}

\vspace{0.5cm} \emph{Acknowledgement.} {\footnotesize The first
author has been partially supported by the Gruppo Nazionale per
l'Analisi Matematica, la Probabilit\`a e le loro Applicazioni
(GNAMPA) of the Istituto Nazionale di Alta Matematica (INdAM) and by
Universit\`a degli Studi di Napoli Parthenope through the project
\lq\lq sostegno alla Ricerca individuale\rq\rq (triennio 2015 - 2017)}.\par

\vspace{4mm}

 \ The second author is grateful to  Yousri Slaoui for sending its a very interest
 article  \cite{Khardani}. \par

\end{document}